\documentclass[12pt,a4paper,reqno]{amsart}
\usepackage{fullpage}
\usepackage{amsmath,amssymb}
\usepackage{braket}
\usepackage{verbatim}
\usepackage{lscape}
\usepackage{newtxtext,newtxmath}
\usepackage{hyperref}
\usepackage{url}
\usepackage{xcolor}
\usepackage{comment}

\newtheorem{teor}{Theorem}
\newtheorem{lem}{Lemma}

\newcommand{\z}{\mathbb}

\newcommand{\dd}{\mathrm d}
\newcommand{\Stilde}{\widetilde{S}}
\newcommand{\M}{\mathfrak{M}}
\newcommand{\m}{\mathfrak{m}}
\newcommand{\gotht}{\mathfrak{t}}
\newcommand{\eps}{\varepsilon}

\allowdisplaybreaks

\author{Alessandro Gambini}

\title{Diophantine approximation with a quaternary problem}

\date{\today}

\allowdisplaybreaks

\begin{document}

\begin{abstract}
Let $1<k<7/6$, $\lambda_1,\lambda_2,\lambda_3$ and $\lambda_4$ be non-zero real numbers, not all of the same sign such that $\lambda_1/\lambda_2$ is irrational and let $\omega$ be a real number. We prove that the inequality $|\lambda_1p_1^2+\lambda_2p_2^2+\lambda_3p_3^2+\lambda_4p_4^k-\omega|\le (\max_j p_j)^{-\frac{7-6k}{14k}+\varepsilon}$ has infinitely many solutions in prime variables $p_1,p_2,p_3,p_4$ for any $\varepsilon>0$.
\end{abstract}

\maketitle

\section{Introduction}

Numerous recent papers have explored a Diophantine inequality involving prime variables, each with a unique set of assumptions and conclusions. In their work, Br\"udern, Cook, and Perelli \cite{brudern-cook-perelli}, focused on binary linear forms in prime arguments. Cook and Fox \cite{cook-fox}, addressed a ternary form with primes squared, and this was subsequently improved in terms of approximation by Harman in \cite{harman2004ternary}. Cook \cite{cook2001general}, provided a more comprehensive description of the problem, which was later refined by Cook and Harman \cite{cook-harman}.

There are several distinctions between the results mentioned above and the scope of our research. Notably, in their papers, the assumption that all coefficients $\lambda_j$ are positive is not a constraint. Additionally, the values of $k_j$ are consistent positive integers for all $j$. However, the pivotal aspect remains the requirement that $\lambda_1/\lambda_2$ must be irrational. In our case, we will prove that there are infinitely many solution to a the problem of the form
\begin{align*}
|\lambda_1p_1^{k_1}+\cdots+\lambda_rp_r^{k_r}-\omega|\le\eta
\end{align*}
when $\eta$ depends on the maximum of the $p_j$, whereas in the previously cited papers, $\eta$ is a small negative power of $\omega$.

Vaughan \cite{vaughan1974diophantineI} follows a similar approach to the one we employ in our article, dealing with a ternary linear form in prime arguments and assuming more suitable conditions on the $\lambda_j$. He proved that there are infinitely many solutions to the problem:

\begin{align*}
|\lambda_1p_1+\lambda_2p_2+\lambda_3p_3-\omega|\le\eta
\end{align*}
when $\eta$ depends on the maximum of the $p_j$. In his case, $\eta=(\max_j p_j)^{-\frac{1}{10}}$. This result was enhanced by Baker and Harman \cite{baker-harman} with an exponent of $-\frac16$, by Harman \cite{harman1991general} with an exponent of $-\frac15$, and finally by Matom\"aki \cite{matomaki} with an exponent of $-\frac29$.

Languasco and Zaccagnini, in \cite{Languasco-Zaccagnini2016} and \cite{languasco-zaccagnini-ternary}, examined a ternary problem with varying powers $k_j$, one of which depended on a parameter $k$. Additionally, Gambini, Languasco, and Zaccagnini \cite{GLZ}, analyzed a ternary problem involving two primes and a $k$-th power of a prime. In all these cases, the value of $\eta$ still depends on the primes $p_j$ also contingent on the parameter $k$. The concept in this scenario is to optimize the value of $k$ to maximize the range in which the inequality holds.

Languasco and Zaccagnini also addressed a quaternary form \cite{languasco-zaccagnini-quaternary} that involved a prime and three squares of primes, resulting in $\eta=(\max_j p_j)^{-\frac{1}{18}}$. This was improved by Li and Wang \cite{li-wang} and later by Liu and Sun in \cite{liu2013diophantine} with $\eta=(\max_j p_j)^{-\frac{1}{16}}$ using the Harman technique. Mu \cite{mu2016diophantine} investigated a problem with five variables comprising four squares of primes and a $k$-th power of a prime, optimizing the value of $k$. Ge and Li \cite{ge-li}, utilized a quaternary form with varying integer powers $k_j$. Gambini  \cite{Gambini}, explored a quaternary problem featuring one prime, two squares of primes, and a $k$-th power of a prime, while Gao and Liu \cite{gao-liu} and later Mu, Zhu, and Li \cite{muzhuli} examined a problem with four squares of a prime and a $k$-th power of a prime.

The case of this paper involves three squares of primes and one $k$-th power of a prime. We prove the following theorem:

\begin{teor}
Assume that $1<k<7/6$, $\lambda_1,\lambda_2,\lambda_3$ and $\lambda_4$ be non-zero real numbers, not all of the same sign, that $\lambda_1/\lambda_2$ is irrational and let $\omega$ be a real number. The inequality
\begin{align}\label{teorema}
\left|\lambda_1p_1^2+\lambda_2p_2^2+\lambda_3p_3^2+\lambda_4p_4^k-\omega\right|\le\left(\max_j p_j\right)^{-\frac{7-6k}{14k}+\varepsilon}
\end{align}
has infinitely many solutions in prime variables $p_1,p_2,p_3,p_4$ for any $\varepsilon>0$.
\end{teor}

\section{Outline of the proof}

We use a variant of the classical circle method that was introduced by Davenport and Heilbronn in 1946 \cite{davenport1946indefinite} substituting the integration over a circle, or equivalently over the interval $[0,1]$, with integration across the whole real line. 

In this paper, we denote prime numbers as $p$ and $p_i$, where $k\geq 1$ is a real number, $\varepsilon$ represents a minute positive value whose specifics might vary depending on occurrences, and $\omega$ is a fixed real number. To establish the existence of infinitely many solutions for \eqref{teorema}, it suffices to create an increasing sequence $X_n$ that grows towards infinity, ensuring that \eqref{teorema} has at least one solution with $\max p_j\in[\delta X_n,X_n]$, with $\delta>0$, a fixed value contingent upon the choice of $\lambda_j$. Consider $q$ as the denominator of a convergent to $\lambda_1/\lambda_2$, with $X_n=X$ (omitting the subscript $n$), and traverses the sequence $X=q^{7/3}$.
We set
\begin{align}
&S_k(\alpha)=\sum_{\delta X\le p^k\le X}\log p\  e(p^k\alpha),\label{S_k}\\
&U_k(\alpha)=\sum_{\delta X\le n^k\le X}e(n^k\alpha), \nonumber\\
&T_k(\alpha)=\int_{(\delta X)^{\frac1k}}^{X^{\frac1k}}e(\alpha t^k)\,\dd t,\label{T_k}
\end{align}
where $e(\alpha)=e^{2\pi i\alpha}$. 

To obtain the most accurate estimate, we utilize the sieve function $\rho(m)$ as defined in (5.2) of \cite{harman-kumchev} introduced by Harman and Kumchev and employed by Wang and Yao in \cite{wang-yao} for the case $k=2$. This function serves as a non-trivial lower bound for the characteristic function of primes. It enables the definition of an exponential function \eqref{stilde} with a distinct weight:
\[
\rho(m)=\psi(m,X^{5/42})-\sum_{X^{5/42}\le p< X^{1/4}}\psi(m/p,z(p)),
\]
where
\[
\psi(m,z)=\left\{\begin{array}{ll}1 & \text{if }p\vert m\ \Rightarrow\ p\ge z, 
\\
0 & \text{otherwise}
\end{array}\right.
\]
and
\[
z(p)=\left\{\begin{array}{ll}X^{5/28}p^{-1/2} & \text{if }p<X^{1/7}, 
\\
p & \text{if } X^{1/7}\le p\le X^{3/14},
\\
X^{5/14}p^{-1} & \text{if }p>X^{3/14}.
\end{array}\right.
\]
The crucial property of $\rho(m)$ we focus on is the estimation (2.3) in \cite{wang-yao}:
\begin{align*}
\sum_{m\in I}\rho(m)=\ell \vert I\vert (\log X)^{-1}+O(X^{1/2}(\log X)^{-2}),
\end{align*}
where $\ell>0$ is an absolute constant and $I$ is any subinterval of $[(\delta X)^{1/2},X^{1/2}]$.
Based on this, we define the following exponential function:
\begin{align}\label{stilde}
\Stilde_2(\alpha)=\sum_{\delta X\le m^2\le X}\rho(m)e(m^2\alpha).
\end{align}

We will approximate $S_k$ with $T_k$ or $U_k$ and we will approximate $\Stilde_2$ with $T_2$.

By the Prime Number Theorem and first derivative estimates for trigonometric integrals we establish
\begin{align}\label{stima_tk}
S_k(\alpha)\ll X^{\frac1k},\qquad \Stilde_2(\alpha)\ll X^{\frac12},\qquad T_k(\alpha)\ll_{k,\delta} X^{\frac1k-1}\min(X,|\alpha|^{-1}),
\end{align}
where $k\ge1$ and $\delta>0$ are real numbers.

Moreover the Euler summation formula implies that, for $k\ge1$,
\begin{align}\label{t-u}
T_k(\alpha)-U_k(\alpha)\ll 1+|\alpha|X.
\end{align}

We also require a continuous function to identify solutions of \eqref{teorema}. Hence, we introduce\begin{align*}
\widehat{K}_\eta(\alpha):=\max\{0,\eta-|\alpha|\}\quad\mbox{where}\quad\eta>0
\end{align*}
whose inverse Fourier transform is
\begin{align*}
K_{\eta}(\alpha)=\left(\frac{\sin(\pi\alpha\eta)}{\pi\alpha}\right)^2
\end{align*}
for $\alpha\neq0$ and, by continuity, $K_{\eta}(0)=\eta^2$. It vanishes at infinity like $|\alpha|^{-2}$ and in fact it is trivial to prove that
\begin{align}\label{k_eta}
K_{\eta}(\alpha)\ll\min(\eta^2,|\alpha|^{-2}).
\end{align}

The original works of Davenport-Heillbronn in \cite{davenport1946indefinite} and later Vaughan in \cite{vaughan1974diophantineI} and \cite{vaughan1974diophantineII} directly approximate the difference $|S_k(\alpha)-T_k(\alpha)|$, estimating it as $O(1)$ using the Euler summation formula. Br\"udern, Cook, and Perelli in \cite{brudern-cook-perelli} enhanced these estimations by computing the $L^2$-norm of $|S_k(\alpha)-T_k(\alpha)|$, leading to substantially improved conditions and a broader major arc compared to the original approach. Introducing the generalized version of the Selberg integral
\begin{align*}
\mathcal{J}_k(X,h)=\int_{X}^{2X}\left(\theta((x+h)^{\frac1k})-\theta(x^{\frac1k})-((x+h)^{\frac1k}-x^{\frac1k})\right)^2\dd x,
\end{align*}
where $\theta$ is the Chebyshev Theta function,
\[
\theta(x)=\sum_{p\le x}\log p,
\]
we have the following lemmas.

\begin{lem}[\cite{Languasco-Zaccagnini2016}, Lemma  1]\label{s-u}
Let $k\ge1$ be a real number. For $0<Y<\frac12$, we have
\begin{align*}
\int_{-Y}^{Y}|S_k(\alpha)-U_k(\alpha)|^2 \dd\alpha\ll_k\frac{{X^{\frac2k-2}}\log^2 X}{Y}+Y^2X+Y^2\mathcal{J}_k\left(X,\frac{1}{2Y}\right).
\end{align*}
\end{lem}
\begin{lem}[\cite{Languasco-Zaccagnini2016}, Lemma  2]\label{jk_stima}
Let $k\ge1$ be a real number and $\varepsilon$ be an arbitrarily small positive constant. There exists a positive constant $c_1(\varepsilon)$, which does not depend on $k$, such that
\begin{align*}
\mathcal{J}_k(X,h)\ll_k h^2X^{\frac2k-1}\exp\left(-c_1\left(\frac{\log X}{\log \log X}\right)^{\frac13}\right)
\end{align*}
uniformly for $X^{1-\frac{5}{6k}+\varepsilon}\le h\le X$.
\end{lem}

\subsection{Setting the problem}
Let 
\begin{align*}\mathcal{P}(X)=\{(p_1,p_2,p_3,p_4): \delta X<p_1^2, p_2^2,p_3^2,p_4^k<X\}
\end{align*} and let us define
\begin{align*}
\mathcal{I}(\eta,\omega,\mathfrak{X})=\int_{\mathfrak{X}}\Stilde_2(\lambda_1\alpha)S_2(\lambda_2\alpha)S_2(\lambda_3\alpha)S_k(\lambda_4\alpha)K_{\eta}(\alpha)e(-\omega\alpha)\dd\alpha
\end{align*}
where $\mathfrak{X}$ is a measurable subset of $\z{R}$.

It follows from the construction of $\rho(m)$ that, if $\omega(m)$ denotes the characteristic function of the set of primes,
\[
\rho(m)\le\omega(m).
\]
Then, from the definitions of $\Stilde_2(\lambda_1\alpha)$ and $S_j(\lambda_i\alpha)$, and performing the Fourier transform for $K_{\eta}(\alpha)$, we obtain
\begin{align*}
\mathcal{I}(\eta,\omega,\z{R})&=\sum_{p_i\in\mathcal{P}(X)}\rho(m_1)\log p_2\log p_3\log p_4\cdot \\&\quad\left(\max(0,\eta-\left|\lambda_1m_1^2+\lambda_2p_2^2+\lambda_3p_3^2+\lambda_4p_4^k-\omega\right|)\right)\nonumber\\ &\leq \eta(\log X)^3\mathcal{N}(X),
\end{align*}
where $\mathcal{N}(X)$ denotes the number of solutions of the inequality \eqref{teorema} with $(p_1,p_2,p_3,p_4)\in\mathcal{P}(X)$. In other words $\mathcal{I}(\eta,\omega,\z{R})$ provides a lower bound for the quantity we are interested in; therefore it is sufficient to prove that $\mathcal{I}(\eta,\omega,\z{R})>0$.

Next, we partition $\mathbb{R}$ into subsets $\M$, $\m$, and $t$, where $\mathbb{R}=\M\cup \m\cup t$, with $\M$ as the major arc, $\m$ as the minor arc (or intermediate arc), and $t$ as the trivial arc, defined as follows:

\begin{align*}
\M=\left[-\frac{P}{X},\frac{P}{X}\right] 
\qquad
\m=\left[-R,-\frac{P}{X}\right]\cup\left[\frac{P}{X},R\right]
\qquad 
t=\z{R}\backslash(\M\cup \m),
\end{align*}
so that  $\mathcal{I}(\eta,\omega,\z{R})=\mathcal{I}(\eta,\omega,\M)+\mathcal{I}(\eta,\omega,\m)+\mathcal{I}(\eta,\omega,t)$.

The parameters $P=P(X)\gg\log (X)$, and $R=R(X)>1/\eta$ are chosen later (see \eqref{cond_major_2} and \eqref{cond_trivial_2}), along with $\eta=\eta(X)$, which, as previously mentioned, we desire to be a small negative power of $\max p_j$ and therefore of $X$ as in \eqref{final_cond2}.

We anticipate having the main term with the correct order of magnitude on $\M$ without any special hypotheses on the coefficients $\lambda_j$. It is crucial to prove that $\mathcal{I}(\eta,\omega,m)$ and $\mathcal{I}(\eta,\omega,t)$ are both $o(\mathcal{I}(\eta,\omega,\M))$: the contribution from the trivial arc is significantly smaller in comparison to the main term. The main challenge lies within the minor arc, where we will require the complete power of the assumptions on the $\lambda_j$ and the theory of continued fractions.

\textbf{Remark}: From this point forward, whenever we use the symbols $\ll$ or $\gg$, we omit the dependence of the approximation on the constants $\lambda_j, \delta$, and $k$.

\subsection{Lemmas}

In this paper we will also use Lemmas 5 of \cite{GLZ} and (2.5) of \cite{wang-yao} that allow us to have an estimation of mean value of $|S_k(\alpha)|^4$ and $|\Stilde_2(\alpha)|^4$:



\begin{lem}[\cite{GLZ}, Lemma 5]\label{s_k^2}
Let $k>1$, $\tau>0$. We have
\begin{align*}
\int_{-\tau}^{\tau}|S_k(\alpha)|^2\dd\alpha
\ll
\left(\tau X^{1/k}+X^{2/k-1}\right)(\log X)^3
\qquad
\int_n^{n+1} |S_k(\alpha)|^2\dd\alpha
\ll
X^{1/k}(\log X)^3.
\end{align*}
\end{lem}


Finally, we will use the following Lemma.
\begin{lem}\label{S_2^4}
\begin{align*}
&\int_0^1 |S_2(\alpha)|^4\dd\alpha\ll X\log^2 X\qquad\int_{\z{R}}  |S_2(\alpha)|^4 K_{\eta}(\alpha)\dd\alpha\ll \eta X\log^2 X.
\\
&\int_0^1 |\Stilde_2(\alpha)|^4\dd\alpha\ll X(\log X)^c\qquad\int_{\z{R}} |\Stilde_2(\alpha)|^4 K_{\eta}(\alpha)\dd\alpha\ll \eta X(\log X)^c.
\end{align*}
\end{lem}

\proof
The first two statements are based on Satz 3 of \cite{rieger}, p. 94 and the last two derive directly from (2.5) of \cite{wang-yao}.
\endproof

\section{The major arc}

We begin with the major arc and the computation of the main term. Substituting all $S_2$, $S_k$, and $\Stilde_2$ defined in \eqref{S_k} and \eqref{stilde} with their respective $T_i$ defined in \eqref{T_k} brings forth some discrepancies that require estimation using Lemma \ref{s-u}, the Cauchy-Schwarz inequality, and the H\"older inequality. We proceed to calculate
\begin{align*}
\mathcal{I}(\eta,\omega,&\M)=\int_{\M}S_2(\lambda_1\alpha)S_2(\lambda_2\alpha)S_2(\lambda_3\alpha)S_k(\lambda_4\alpha)K_{\eta}(\alpha)e(-\omega\alpha)\dd\alpha\nonumber\\
=&\ell(\log X)^{-1}\int_{\M}T_2(\lambda_1\alpha)T_2(\lambda_2\alpha)T_2(\lambda_3\alpha)T_k(\lambda_4\alpha)K_{\eta}(\alpha)e(-\omega\alpha)\dd\alpha\nonumber\\
&+\int_{\M}(\Stilde_2(\lambda_1\alpha)-\ell(\log X)^{-1}T_2(\lambda_1\alpha))S_2(\lambda_2\alpha)T_2(\lambda_3\alpha)T_k(\lambda_4\alpha)K_{\eta}(\alpha)e(-\omega\alpha)\dd\alpha\nonumber\\
&+\int_{\M}\Stilde_2(\lambda_1\alpha)(S_2(\lambda_2\alpha)-T_2(\lambda_2\alpha))T_2(\lambda_3\alpha)T_k(\lambda_4\alpha)K_{\eta}(\alpha)e(-\omega\alpha)\dd\alpha\nonumber\\
&+\int_{\M}\Stilde_2(\lambda_1\alpha))S_2(\lambda_2\alpha)(S_2(\lambda_3\alpha)-T_2(\lambda_3\alpha))T_k(\lambda_4\alpha)K_{\eta}(\alpha)e(-\omega\alpha)\dd\alpha\nonumber\\
&+\int_{\M}\Stilde_2(\lambda_1\alpha)S_2(\lambda_2\alpha)S_2(\lambda_3\alpha)(S_k(\lambda_4\alpha)-T_k(\lambda_4\alpha))K_{\eta}(\alpha)e(-\omega\alpha)\dd\alpha\nonumber\\
=&J_1+J_2+J_3+J_4+J_5,
\end{align*}
say. Since the computations for $J_3$ is similar to, but simpler than, the corresponding ones for $J_2$, $J_4$ and $J_5$, we will leave it to the reader.

\subsection{Main Term: lower bound for $J_1$}

As the reader might expect the main term is given by the summand $J_1$.

Let $H(\alpha)=T_1(\lambda_1\alpha)T_2(\lambda_2\alpha)T_2(\lambda_3\alpha)T_k(\lambda_4\alpha)K_{\eta}(\alpha)e(-\omega\alpha)$ so that

\begin{align*}
J_1=\ell(\log X)^{-1}\int_{\z{R}}H(\alpha)\dd\alpha+\mathcal{O}\left(\int_{P/X}^{+\infty}|H(\alpha)|\dd\alpha\right).
\end{align*}

Using inequalities \eqref{stima_tk} and \eqref{k_eta} ,

\begin{align*}
\int_{P/X}^{+\infty}|H(\alpha)|\dd\alpha\ll& X^{-\frac12}X^{-\frac12} X^{-\frac12} X^{\frac1k-1}\eta^2\int_{P/X}^{+\infty}\frac{\dd\alpha}{\alpha^4}\nonumber\\\ll& X^{\frac1k-\frac52}\eta^2\frac{X^3}{P^3}=X^{\frac1k+\frac12}\eta^2P^{-3}=o\left(X^{\frac1k+\frac12}\eta^2\right)
\end{align*}
provided that $P\rightarrow+\infty$. Let $D=[(\delta X)^{\frac12},X^{\frac12}]^3\times[(\delta X)^{\frac1k},X^{\frac1k}]$ we have
\begin{align*}
\int_{\z{R}}H(\alpha)\dd\alpha&=\idotsint_D\int_{\z{R}}e((\lambda_1t_1^2+\lambda_2t_2^2+\lambda_3t_3^2+\lambda_4t_4^k-\omega)\alpha)K_{\eta}(\alpha)\dd\alpha\,\dd t_1\dd t_2\dd t_3\dd t_4\\&=\idotsint_D\max(0,\eta-|\lambda_1t_1^2+\lambda_2t_2^2+\lambda_3t_3^2+\lambda_4t_4^k-\omega)|)\dd t_1\dd t_2\dd t_3\dd t_4.
\end{align*}

Apart from very slight changes in the computation, we proceed with a change of variables as in \cite{Gambini} and we obtain 
\begin{align*}
J_1
\gg
(\log X)^{-1}\eta^2X^{\frac1k+\frac12},
\end{align*}
which is the expected lower bound.

\subsection{Bound for \texorpdfstring{$J_2$}{J2}}

We expect the main term to have the dominant asymptotic behavior, then we shall prove that all the remaining terms of the sum are $o\left((\log X)^{-1}\eta^2 X^{\frac1k+\frac12}\right)$.

From partial summation on \eqref{stilde} we get
\[
\Stilde_2(\lambda_1\alpha)=\int_{(\delta X)^{\frac12}}^{X^{\frac12}}e(\lambda t^2\alpha)\,\mathrm{d}\bigg(\sum_{\substack{m_2\le t\\ m_2\in[(\delta X)^{1/2},X^{1/2}]}}\rho(m_2)\bigg),
\]
then
\[
\Stilde_2(\lambda_1\alpha)-\ell (\log X)^{-1}T_2(\lambda_1 \alpha)\ll X^{\frac12}(\log X)^{-2}(1+|\alpha|X).
\]

Retrieving \eqref{k_eta} and using the triangle inequality,
\begin{align*}
J_2\ll
&\eta^2\int_{\M}|\Stilde_2(\lambda_1\alpha)-\ell (\log X)^{-1}T_2(\lambda_1\alpha)||T_2(\lambda_2\alpha)||T_2(\lambda_3\alpha)||T_k(\lambda_4\alpha)|\dd\alpha\nonumber\\
\ll&\eta^2X^{\frac12}(\log X)^{-2}(1+|\alpha|X)\int_{\M}|T_2(\lambda_2\alpha)||T_2(\lambda_3\alpha)||T_k(\lambda_4\alpha)|\dd\alpha
\\
\ll& \eta^2X^{\frac12}(\log X)^{-2}\int_{0}^{1/X}|T_2(\lambda_1\alpha)||T_2(\lambda_3\alpha)||T_k(\lambda_4\alpha)|\dd\alpha
\\
&+\eta^2X^{\frac32}(\log X)^{-2}\int_{1/X}^{P/X}\alpha|T_2(\lambda_1\alpha)||T_2(\lambda_3\alpha)||T_k(\lambda_4\alpha)|\dd\alpha
\\
=&o\left(\eta^2X^{\frac1k+\frac12}(\log X)^{-1}\right).
\end{align*}

\subsection{Bound for \texorpdfstring{$J_4$}{J4}}

Using the triangle inequality and \eqref{k_eta},
\begin{align*}
J_4=&\int_{\M}\Stilde_2(\lambda_1\alpha)S_2(\lambda_2\alpha)(S_2(\lambda_3\alpha)-T_2(\lambda_3\alpha))T_k(\lambda_4\alpha)K_{\eta}(\alpha)e(-\omega\alpha)\dd\alpha\nonumber\\
\ll&\eta^2\int_{\M}|\Stilde_2(\lambda_1\alpha)||S_2(\lambda_2\alpha)||S_2(\lambda_3\alpha)-T_2(\lambda_3\alpha)||T_k(\lambda_4\alpha)|\dd\alpha\nonumber\\ 
\le&\eta^2\int_{\M}|\Stilde_2(\lambda_1\alpha)||S_2(\lambda_2\alpha)||S_2(\lambda_3\alpha)-U_2(\lambda_3\alpha)||T_k(\lambda_4\alpha)|\dd\alpha\nonumber\\
&+\eta^2\int_{\M}|\Stilde_2(\lambda_1\alpha)||S_2(\lambda_2\alpha)||U_2(\lambda_3\alpha)-T_2(\lambda_3\alpha)||T_k(\lambda_4\alpha)|\dd\alpha\nonumber\\
=&\eta^2(A_4+B_4),
\end{align*}
say. Using Theorem \ref{stima_tk} and the H\"older inequality,
\begin{align*}
A_4\ll &X^{\frac1k}\int_{\M}|\Stilde_2(\lambda_1\alpha)||S_2(\lambda_2\alpha)||S_2(\lambda_3\alpha)-U_2(\lambda_3\alpha)|\dd\alpha\nonumber\\
\ll & X^{\frac1k}\left(\int_{\M}|\Stilde_2(\lambda_1\alpha)|^4\dd\alpha\right)^{\frac14}\left(\int_{\M}|S_2(\lambda_2\alpha)|^4\dd\alpha\right)^{\frac14}\left(\int_{\M}|S_2(\lambda_3\alpha)-U_2(\lambda_3\alpha)|^2\dd\alpha\right)^{\frac12}.
\end{align*}

Using Lemmas \ref{s-u}-\ref{jk_stima}-\ref{S_2^4}, for any fixed $A$,
\begin{align*}
A_4\ll X^{\frac1k}(X\log^2 X)^{\frac12}(\log X)^{-\frac{A}{2}}=X^{\frac12+\frac1k}(\log X)^{1-\frac{A}{2}}=o\left((\log X)^{-1}X^{\frac1k+\frac12}\right)
\end{align*}
as long as $A>4$.

As for $A_2$ we used in the estimation above Lemma \ref{s-u} that has two more terms, but also in this case these terms are negligible if we want to meet the hypothesis of Lemma \ref{jk_stima}: in fact it requires that 
\[
X^{1-\frac5{12}+\varepsilon} \le \frac{X}{P} \leq X
\]
and this is consistent with the choice we will make in \eqref{cond2_major_2}.

Again using Theorem \ref{t-u},
\begin{align*}
B_4=&\int_{\M}|\Stilde_2(\lambda_1\alpha)||S_2(\lambda_2\alpha)||U_2(\lambda_3\alpha)-T_2(\lambda_3\alpha)||T_k(\lambda_4\alpha)|\dd\alpha\nonumber\\
\ll &\int_0^{1/X}|\Stilde_2(\lambda_1\alpha)||S_2(\lambda_2\alpha)||T_k(\lambda_4\alpha)|\dd\alpha\nonumber\\
&+X\int_{1/X}^{P/X}\alpha|\Stilde_2(\lambda_1\alpha)||S_2(\lambda_2\alpha)||T_k(\lambda_4\alpha)|\dd\alpha.
\end{align*}

Remembering that $|\alpha|\le\frac{P}{X}$ on $\M$ and using the H\"older inequality, trivial bounds and Lemma \ref{S_2^4} we have
\begin{align*}
B_4\ll &X^{\frac12}\,X^{\frac12}X^{\frac1k}\frac1X+X\,X^{\frac1k}\left(\int_{1/X}^{P/X}\alpha^2\right)^{\frac12}\left(\int_{1/X}^{P/X}|\Stilde_2(\lambda_1\alpha)|^4\dd\alpha\right)^{\frac14}\left(\int_{1/X}^{P/X}|S_2(\lambda_2\alpha)|^4\dd\alpha\right)^{\frac14}\nonumber\\
\ll & X^{\frac1k}+X^{1+\frac1k}\left(X\log^2 X\right)^{\frac12}\left(\int_{1/X}^{P/X}\alpha^2\dd\alpha\right)^{\frac12}\nonumber\\
\ll & X^{\frac1k}+X^{\frac32+\frac1k}\log X\left(\frac{P}{X}\right)^{\frac32}=X^{\frac1k}P^{\frac32}\log X.
\end{align*}

We assume
\begin{align}\label{cond2_major_2}
P\le X^{\frac13-\varepsilon},
\end{align}
so that
$P^{\frac32}=o(X^{\frac12}/\log^2 X)$ which, with the upper bound for $B_4$ here above, ensures that 
\[
B_4=o((\log X)^{-1}X^{1/2+1/k}).
\]

\subsection{Bound for \texorpdfstring{$J_5$}{J5}}

In order to provide an estimation for $J_5$, we use \eqref{k_eta},
\begin{align*}
J_5\ll&\eta^2\int_{\M}|\Stilde_2(\lambda_1\alpha)||S_2(\lambda_2\alpha)||S_2(\lambda_3\alpha)||S_k(\lambda_4\alpha)-T_k(\lambda_4\alpha)|\dd\alpha
\end{align*}
and then the arithmetic-geometric inequality:
\begin{align*}
J_5\ll&\eta^2\sum_{j=2}^3\left(\int_{\M}|\Stilde(\lambda_1\alpha)||S_2(\lambda_j\alpha)|^2|S_k(\lambda_4\alpha)-T_k(\lambda_4\alpha)|\dd\alpha\right).
\end{align*}

The three terms may be estimated in the same way and produce the same
upper bound. We show the details of the bound only for the case $j=2$:
\begin{align*}
\eta^2\int_{\M}&|\Stilde(\lambda_1\alpha)||S_2(\lambda_2\alpha)|^2|S_k(\lambda_4\alpha)-T_k(\lambda_4\alpha)|\dd\alpha\nonumber\\
\ll &\eta^2\int_{\M}|\Stilde(\lambda_1\alpha)||S_2(\lambda_2\alpha)|^2|S_k(\lambda_4\alpha)-U_k(\lambda_4\alpha)|\dd\alpha\nonumber\\ 
&+\eta^2\int_{\M}|\Stilde(\lambda_1\alpha)||S_2(\lambda_2\alpha)|^2|U_k(\lambda_4\alpha)-T_k(\lambda_4\alpha)|\dd\alpha\nonumber\\ 
=&\eta^2(A_5+B_5),
\end{align*}
say. Using trivial estimates,
\begin{align*}
A_5\ll X^{\frac12}\int_{\M}|S_2(\lambda_2\alpha)|^2|S_k(\lambda_4\alpha)-U_k(\lambda_4\alpha)|\dd\alpha
\end{align*}
then using the Cauchy-Schwartz inequality, for any fixed $A>4$, by Lemmas \ref{S_2^4}, \ref{s-u} and \ref{jk_stima} we have
\begin{align*}
A_5\ll &X^{\frac12}\left(\int_{\M}|S_2(\lambda_1\alpha)|^4\dd\alpha\right)^{\frac12}\left(\int_{\M}|S_k(\lambda_4\alpha)-U_k(\lambda_4\alpha)|^2\dd\alpha\right)^{\frac12}\nonumber\\\ll& X^{\frac12}\, X^{\frac12}\log X \frac{P}{X}\mathcal{J}_k\left(X,\frac{X}{P}\right)^{\frac12}\ll_A X^{\frac12+\frac1k}(\log X)^{1-\frac{A}{2}}
=o\left((\log X)^{-1}X^{\frac12+\frac1k}\right)
\end{align*}
provided that $\frac{X}{P}\ge X^{1-\frac{5}{6k}+\varepsilon}$ (condition of Lemma \ref{jk_stima}), that is, 
\begin{align}\label{cond3_major_2}
(\log X)^A\ll_A P\le X^{\frac{5}{6k}-\varepsilon}.
\end{align}

Now we turn to $B_5$, using Theorem \ref{t-u}:
\begin{align*}
B_5\ll \int_0^{1/X}|\Stilde(\lambda_1\alpha)||S_2(\lambda_2\alpha)|^2\dd\alpha+X\int_{1/X}^{P/X}\alpha|\Stilde(\lambda_1\alpha)||S_2(\lambda_2\alpha)|^2\dd\alpha.
\end{align*}

Using trivial estimates and Lemma \ref{S_2^4},
\begin{align*}
B_5\ll&X^{\frac32}\frac1X+X\cdot X^{\frac12}\left(\int_{1/X}^{P/X}\alpha^2\dd\alpha\cdot\int_{1/X}^{P/X}|S_2(\lambda_1\alpha)|^4\dd\alpha\right)^{\frac12}\nonumber\\
\ll& X^{\frac12}+X^{\frac32}(P/X)^{\frac32}X^{\frac12}\log X=X^{\frac12}+P^{\frac32}X^{\frac12}\log X.
\end{align*}

The case $j=3$ can be estimated in the same way. We need
\begin{align*}
P=o\left(X^{\frac{2}{3k}-\varepsilon}\right).
\end{align*}

Summing up with \eqref{cond3_major_2}, 
\begin{align}\label{cond4_major_2}
P\le X^{\frac{2}{3k}-\varepsilon}.
\end{align}

Collecting all the bounds for $P$, that is, \eqref{cond2_major_2}, \eqref{cond3_major_2}, \eqref{cond4_major_2} we can take
\begin{align}\label{cond_major_2}
P\le X^{\frac13-\varepsilon}.
\end{align}

In fact, if we consider \eqref{cond2_major_2}, \eqref{cond3_major_2} and \eqref{cond4_major_2} we should choose the most restrictive condition among the three but as we expect that the value of $k$ is smaller than 2, \eqref{cond2_major_2} is the most restrictive: $\frac{2}{3k}\le\frac{5}{6k}$ and $\frac13\ge\frac{2}{3k}$ only if $k\ge2$.

\section{Trivial arc}

By the arithmetic-geometric mean inequality and the trivial bound for $\Stilde_2(\lambda_1\alpha)$, we see that
\begin{align*}
|\mathcal{I}(\eta,\omega,\gotht)|\ll&\int_{R}^{+\infty}|\Stilde_2(\lambda_1\alpha)S_2(\lambda_2\alpha)S_2(\lambda_3\alpha)S_k(\lambda_4\alpha)K_{\eta}(\alpha)|\dd\alpha\nonumber\\\ll&X^{\frac12}\sum_{j=2}^3\int_{R}^{+\infty}|S_2(\lambda_j\alpha)|^2|S_k(\lambda_4\alpha)|K_{\eta}(\alpha)\dd\alpha.
\end{align*}

The three terms may be estimated in the same way and produce the same upper bound. We show the details of the bound only for the case $j=1$:
\begin{align*}
X^{\frac12}\int_{R}^{+\infty}&|S_2(\lambda_j\alpha)|^2|S_k(\lambda_4\alpha)|K_{\eta}(\alpha)\dd\alpha\\
\ll& X^{\frac12}\left(\int_{R}^{+\infty}|S_2(\lambda_1\alpha)|^4K_{\eta}(\alpha)\dd\alpha\right)^{\frac12}\left(\int_{R}^{+\infty}|S_k(\lambda_4\alpha)|^2K_{\eta}(\alpha)\dd\alpha\right)^{\frac12}\\
\ll& X^{\frac12}\left(\int_{R}^{+\infty}\frac{|S_2(\lambda_1\alpha)|^4}{\alpha^2}\dd\alpha\right)^{\frac12}\left(\int_{R}^{+\infty}\frac{|S_k(\lambda_4\alpha)|^2}{\alpha^2}\dd\alpha\right)^{\frac12}=X^{\frac12}C_1^{\frac12}C_2^{\frac12},
\end{align*}
say. Using Lemma \ref{S_2^4}, we have
\begin{align}\label{c1_2}
C_1=&\int_{R}^{+\infty}\frac{|S_2(\lambda_2\alpha)|^4}{\alpha^2}\dd\alpha\ll\int_{\lambda_1R}^{+\infty}\frac{|S_2(\alpha)|^4}{\alpha^2}\dd\alpha\nonumber\\\ll&\sum_{n\ge|\lambda_1| R}\frac{1}{(n-1)^2}\int_{n-1}^n |S_2(\alpha)|^4\dd\alpha\ll\frac{X\log^2 X}{R}.
\end{align}

Now using Lemma \ref{s_k^2},
\begin{align}\label{c2_2}
C_2=&\int_{R}^{+\infty}\frac{|S_k(\lambda_4\alpha)|^4}{\alpha^2}\dd\alpha\ll\int_{\lambda_4R}^{+\infty}\frac{|S_k(\alpha)|^2}{\alpha^2}\dd\alpha\nonumber\\\ll&\sum_{n\ge|\lambda_4| R}\frac{1}{(n-1)^2}\int_{n-1}^n |S_k(\alpha)|^2\dd\alpha\ll\frac{X^{\frac1k}\log^3 X}{R}.
\end{align}

Collecting \eqref{c1_2} and \eqref{c2_2},
\begin{align*}
|\mathcal{I}(\eta,\omega,\gotht)|\ll X^{\frac12}\left(\frac{X\log^2 X}{R}\right)^{\frac12}\left(\frac{X^{\frac1k}\log^3 X}{R}\right)^{\frac12}\ll\frac{X^{1+\frac1{2k}}(\log X)^{\frac52}}{R}.
\end{align*}

Hence, remembering that $|\mathcal{I}(\eta,\omega,t)|$ must be $o\left((\log X)^{-1}\eta^2X^{\frac1k+1}\right)$, i.e. little-o of the main term, the choice
\begin{align}\label{cond_trivial_2}
R=\frac{X^{\frac12-\frac1{2k}}\log^4 X}{\eta^2}
\end{align}
is admissible.

\section{The minor arc}

In \cite{wang-yao} section 4 it is proven that the measure of the set where $|\Stilde_2(\lambda_1\alpha)|^{\frac12}$ and $|\Stilde_2(\lambda_2\alpha)|$ are both large for $\alpha\in m$ is small, exploiting the fact that the ratio $\lambda_1/\lambda_2$ is irrational.

\begin{lem}[Wang-Yao \cite{wang-yao}, Lemma 1]\label{harman_s2}
Suppose that $X^{\frac12}\ge Z\ge X^{\frac12-\frac{1}{14}+\varepsilon}$ and $|\Stilde_2(\lambda\alpha)|>Z$. Then there are coprime integers $(a,q)=1$ satisfying
\begin{align*}
1\le q\le \left(\frac{X^{\frac12+\varepsilon}}{Z}\right)^4,\qquad |q\lambda\alpha-a|\ll X^{-1}\left(\frac{X^{\frac12+\varepsilon}}{Z}\right)^4.
\end{align*}
\end{lem}

In this case we need only Lemma \ref{harman_s2}. Let us split $\m$ into subsets $\m_1$, $\m_2$ and $\m^*=\m\backslash(\m_1\cup\m_2)$ where
\begin{align*}
&\m_i=\{\alpha\in \m\,:\ |\Stilde_2(\lambda_i\alpha)|\le X^{\frac12-u+\varepsilon}\} 
\end{align*}
remembering that Lemma \ref{harman_s2} holds for $0\le u\le \frac{1}{14}$. In this case we leave only the parameter $u$ free.
Using the H\"older inequalities and the definition of $m_i$ we obtain
\begin{align}\label{cond_m1_2}
\left|\mathcal{I}(\eta,\omega,\m_i)\right|\ll&\int_{\m_i}|\Stilde_2(\lambda_1\alpha)||S_2(\lambda_2\alpha)||S_2(\lambda_3\alpha)||S_k(\lambda_4\alpha)|K_{\eta}(\alpha)\dd\alpha\nonumber\\ 
\ll&\max{|\Stilde_2(\lambda_1\alpha)|}\left(\int_{\m_i}|S_2(\lambda_2\alpha)|^4K_{\eta}(\alpha)\dd\alpha)\right)^{1/4}\nonumber\\
&\left(\int_{\m_i}|S_2(\lambda_3\alpha)|^4K_{\eta}(\alpha)\dd\alpha)\right)^{1/4}\left(\int_{\m_i}|S_k(\lambda_4\alpha)|^2K_{\eta}(\alpha)\dd\alpha)\right)^{1/2}\nonumber\\
\ll& X^{\frac12-u+\varepsilon}(\eta X\log^2 X)^{\frac14}(\eta X\log^2 X)^{\frac14}\left(\eta X^{\frac1k}\log^3 X\right)^{\frac12}\nonumber\\
=&\eta X^{1-u+\frac{1}{2k}+\varepsilon}\log^{\frac52}X.
\end{align}
by Lemma \ref{harman_s2}.
The bound \eqref{cond_m1_2} must be $o\left(\eta^2X^{\frac12+\frac1k}\right)$, consequently we have the following condition:
\begin{align*}
\eta=\infty\left(X^{{\frac12}-\frac{1}{2k}-u+\varepsilon}\right),
\end{align*}
where we used the notation $f=\infty(g)$ for $g=o(f)$.

It remains to discuss the set $\m^*$ in which the following bounds hold simultaneously
\begin{align*}
|\Stilde_2(\lambda_i\alpha)|>X^{\frac12-u+\varepsilon},\qquad \frac{P}{X}=X^{-\frac23}<|\alpha|\le \frac{\log^2 X}{\eta^2}.
\end{align*} 

Following the dyadic dissection argument shown by Harman in \cite{harman2004ternary} we divide $\m^*$ into disjoint sets $E(Z_1,Z_2,y)$ in which, for $\alpha\in E(Z_1,Z_2,y)$, we have
\begin{align*}
Z_1<|\Stilde_2(\lambda_1\alpha)|\le 2Z_1,\qquad Z_2<|\Stilde_2(\lambda_2\alpha)|\le 2Z_2,\qquad y<|\alpha|\le 2y
\end{align*}
where $Z_i=2^{k_i}X^{\frac12-u+\varepsilon}$ for  $i=1,2$, and $y=2^{k_3}X^{-\frac23-\varepsilon}$ for some non-negative integers $k_1,k_2,k_3$.

It follows that the disjoint sets are, at most, $\ll\log^3 X$. Let us define $\mathcal{A}$ a shorthand for the sets $E(Z_1,Z_2,y)$; we have the following result about the Lebesgue measure of $\mathcal{A}$ following the same lines of Lemma 6 in \cite{mu2016diophantine}. 

In the subsequent Lemma, it is essential for both integers $a_1$ and $a_2$ involved in \eqref{misura2} below not to equal zero: specifically, if $a_1=0$, for instance, then $q_1=1$, and $\vert\alpha\vert$ becomes so small that it cannot belong to $\m^*$. Consequently, to apply the Harman technique, we are compelled to move away from the major arc, where $a_1a_2=0$. As we shall later observe, upon defining the parameter $u$, we won't encounter a gap between the major and minor arcs, obviating the necessity to introduce an intermediate arc.

\begin{lem}\label{lemma_mu2}
We have
\begin{align*}
\mu(\mathcal{A})\ll yX^{2+8u+3\varepsilon}Z_1^{-4}Z_2^{-4}
\end{align*}
where $\mu(\cdot)$ denotes the Lebesgue measure.
\end{lem}
\proof
If $\alpha\in\mathcal{A}$, by Lemma \ref{harman_s2} there are coprime integers $(a_1,q_1)$ and $(a_2,q_2)$ such that
\begin{align}\label{misura2}
1\le q_2\ll\left(\frac{X^{\frac12+\varepsilon/4}}{Z_2}\right)^4,\qquad |q_2\lambda_2\alpha-a_2|\ll X^{-1}\left(\frac{X^{\frac12+\varepsilon/4}}{Z_2}\right)^4
\end{align}

We remark that $a_1a_2\neq0$ otherwise we would have $\alpha\in\M$. In fact, if $a_i=0$ recalling the definitions of $Z_1$ and \eqref{misura2}, we get
\begin{align*}
\vert\alpha\vert\ll q_2^{-1}X^{-1}\left(\frac{X^{\frac12+\varepsilon/4}}{Z_2}\right)^4\ll \frac{X}{X^{2-4u+3\varepsilon}}=X^{-1+4u-3\varepsilon}.
\end{align*}

It means that, on the minor arc
\[
\vert\alpha\vert\gg X^{-1+4u-3\varepsilon}.
\]

We wonder now if there is a gap between the end of the major arc and the beginning of the minor arc: from Lemma \ref{harman_s2} we are sure that $u\le\frac{1}{16}$; furthermore, from the previous lower bound for $\alpha$, we need to check whether $\frac{P}{X}$ is greater than it:
\begin{align*}
&X^{-1+4u-3\varepsilon}<\frac{P}{X}=X^{-\frac23}\quad\Rightarrow\quad u<\frac{1}{12}.
\end{align*}

It is clear that we can choose any parameter $u$ with the condition given by Lemma \ref{harman_s2} without leaving any gap from the two arcs. 

Now, we can further split $\m^*$ into sets $I(Z_1,Z_2,y,Q_1,Q_2)$ where \mbox{$Q_j\le q_j\le 2Q_j$} on each set. Note that $a_i$ and $q_i$ are uniquely determined by $\alpha$. In the opposite direction, for a given quadruple $a_1$, $q_1$, $a_2$, $q_2$ the inequalities
\eqref{misura2} define an interval of $\alpha$ of length
\begin{align*}
\mu(I)\ll\min\left(Q_1^{-1}X^{-1}\left(\frac{X^{\frac12+\varepsilon/4}}{Z_1}\right)^4,Q_2^{-1}X^{-1}\left(\frac{X^{\frac12+\varepsilon/4}}{Z_2}\right)^4\right).
\end{align*}

Taking the geometric mean ($\min(a,b)\le\sqrt{a}\sqrt{b}$) we can write
\begin{align}\label{lower_Q2}
\mu(I)\ll Q_1^{-\frac12}Q_2^{-\frac12}X^{-1}\left(\frac{X^{\frac12+\varepsilon/4}}{Z_1}\right)^2\left(\frac{X^{\frac12+\varepsilon/4}}{Z_2}\right)^2\ll \frac{X^{1+\varepsilon}}{Q_1^{\frac12}Q_2^{\frac12}Z_1^2Z_2^2}.
\end{align}

Now we need a lower bound for $Q_1^{\frac12}Q_2^{\frac12}$: by \eqref{misura2}
\begin{align*}
\left|a_2q_1\frac{\lambda_1}{\lambda_2}-a_1q_2\right|&=\left|\frac{a_2}{\lambda_2\alpha}(q_1\lambda_1\alpha-a_1)-\frac{a_1}{\lambda_2\alpha}(q_2\lambda_2\alpha-a_2)\right|\nonumber\\&\ll q_2|q_1\lambda_1\alpha-a_1|+q_1|q_2\lambda_2\alpha-a_2|\nonumber\\&\ll Q_2X^{-1}\left(\frac{X^{\frac12+\varepsilon/4}}{Z_1}\right)^4+Q_1X^{-1}\left(\frac{X^{\frac12+\varepsilon/4}}{Z_2}\right)^4.
\end{align*}

Remembering that $Q_i\ll \left(\frac{X^{\frac12+\varepsilon/4}}{Z_i}\right)^4$, $Z_i\gg X^{\frac12-u+\varepsilon}$, 
\begin{align}\label{harman_stima2}
\left|a_2q_1\frac{\lambda_1}{\lambda_2}-a_1q_2\right|\ll\left(\frac{X^{\frac12+\varepsilon/4}}{X^{\frac12-u+\varepsilon}}\right)^4X^{-1}\left(\frac{X^{\frac12+\varepsilon/4}}{X^{\frac12-u+\varepsilon}}\right)^4\ll\frac{X^{3+2\varepsilon}}{X^{4-8u+8\varepsilon}}\ll X^{-1+8u-6\varepsilon}.
\end{align}

We recall that $q = X^{1-8u}$ is a denominator of a convergent of $\lambda_1/\lambda_2$. Hence by \eqref{harman_stima2}
Legendre's law of best approximation implies that $\vert a_2q_1\vert \ge q$ and by the same token, for any pair $\alpha$, $\alpha'$ having distinct
associated products $a_2q_1$ (see \cite{watson}, Lemma 2),
\[
\vert a_2(\alpha)q_1(\alpha)-a_2(\alpha')q_1(\alpha')\vert \ge q;
\]
thus, by the pigeon-hole principle, there is at most one value of $a_2q_1$ in the 
interval $[rq,(r+1)q)$ for any positive integer $r$. Hence $a_2q_1$ determines $a_2$ and 
$q_1$ to within $X^{\eps}$ possibilities (the upper bound for the divisor function) and consequently also $a_2q_1$ determines $a_1$ and $q_2$ to within $X^{\eps}$ 
 possibilities from \eqref{harman_stima2}.

Hence we got a lower bound for $q_1q_2$, remembering that in our shorthand $Q_j\le q_j\le 2Q_j$:
\begin{align*}
q_1q_2=a_2q_1\frac{q_2}{a_2}\gg\frac{rq}{|\alpha|}\gg rqy^{-1}
\end{align*}
for the quadruple under consideration.  As a consequence we obtain from \eqref{lower_Q2}, that the total length of the interval $I(Z_1,Z_2,y,Q_1,Q_2)$ with $a_2q_1\in[rq,(r+1)q)$ does not exceed
\begin{align*}
\mu(I)\ll X^{1+2\varepsilon}{Z_1^{-2}Z_2^{-2}}r^{-\frac12}q^{-\frac12}y^{\frac12}.
\end{align*}

Now we need a bound for $r$: inside the interval $[rq,(r+1)q)$, $rq\le |a_2q_1|$ and, in turn from \eqref{misura2}, $a_2\ll q_2|\alpha|$, then
\begin{align*}
&rq\ll q_1q_2|\alpha|\ll\left(\frac{X^{\frac12+\varepsilon/4}}{Z_1}\right)^4\left(\frac{X^{\frac12+\varepsilon/4}}{Z_2}\right)^4y\ll yX^{4+2\varepsilon}Z_1^{-4}Z_2^{-4}\nonumber\\&\Rightarrow r\ll q^{-1}yX^{4+2\varepsilon}Z_1^{-4}Z_2^{-4}.
\end{align*}

Now, we sum on every interval to get an upper bound for the measure of $\mathcal{A}$:
\begin{align*}
\mu(\mathcal{A})\ll X^{1+2\varepsilon}{Z_1^{-2}Z_2^{-2}}q^{-\frac12}y^{\frac12}\sum_{1\le r\ll q^{-1}yX^{4+2\varepsilon}Z_1^{-4}Z_2^{-4}}r^{-\frac12}.
\end{align*}

By standard estimation we obtain
\begin{align*}
\sum_{1\le r\ll q^{-1}yX^{4+2\varepsilon}Z_1^{-4}Z_2^{-4}}r^{-\frac12}\ll (q^{-1}yX^{4+2\varepsilon}Z_1^{-4}Z_2^{-4})^{\frac12}
\end{align*}
then
\begin{align*}
\mu(\mathcal{A})\ll yX^{3+3\varepsilon}Z_1^{-4}Z_2^{-4}q^{-1}\ll yX^{3+3\varepsilon}Z_1^{-4}Z_2^{-4}X^{-1+8u}\ll yX^{2+8u+3\varepsilon}Z_1^{-4}Z_2^{-4}.
\end{align*}

This concludes the proof of Lemma \ref{lemma_mu2}.
\endproof

Using Lemma \ref{lemma_mu2} we finally are able to get a bound for $\mathcal{I}(\eta,\omega,\mathcal{A})$:
\begin{align*}
\mathcal{I}(\eta,\omega,\mathcal{A})=&\int_{\m^*}|\Stilde_2(\lambda_1\alpha)||\Stilde_2(\lambda_2\alpha)||S(\lambda_3\alpha)||S_k(\lambda_4\alpha)|K_{\eta}(\alpha)\dd\alpha\nonumber\\
\ll&\left(\int_{\mathcal{A}}|\Stilde_2(\lambda_1\alpha)\Stilde_2(\lambda_2\alpha)|^4K_{\eta}(\alpha)\dd\alpha\right)^{\frac14}\left(\int_{\mathcal{A}}|S_2(\lambda_3\alpha)|^4K_{\eta}(\alpha)\dd\alpha\right)^{\frac14}\nonumber\\
&\left(\int_{\mathcal{A}}|S_k(\lambda_4\alpha)|^2K_{\eta}(\alpha)\dd\alpha\right)^{\frac12}\nonumber\\
\ll&\left(\min\left(\eta^2,\frac{1}{y^2}\right)\right)^{\frac14}\left((Z_1Z_2)^4\mu(\mathcal{A})\right)^{\frac14}\left(\eta X \log^2 X\right)^{\frac14}\left(\eta X^{\frac1k}\log^3 X\right)^{\frac12}\nonumber\\
\ll&\left(\min\left(\eta^2,\frac{1}{y^2}\right)\right)^{\frac14}Z_1Z_2(yX^{2+8u+4\varepsilon}Z_1^{-4}Z_2^{-4})^{\frac14}\eta^{\frac34}X^{\frac14+\frac{1}{2k}+\varepsilon}\nonumber\\
\ll&\left(\min\left(\eta^2,\frac{1}{y^2}\right)\right)^{\frac14}y^{\frac14}\eta^{\frac34}X^{\frac34+2u+\frac{1}{2k}+\varepsilon}\nonumber\\
\ll&\eta X^{\frac34+u+\frac{1}{2k}+\varepsilon}
\end{align*}
and this must be $o\left(X^{{1+\frac1k}-\varepsilon}\right)$.

The condition on $\eta$ is
\begin{align}\label{final_cond2}
\eta=\infty\left(X^{{\frac14}-\frac{1}{2k}+2u+\varepsilon}\right).
\end{align}

Collecting all the conditions \eqref{cond_m1_2}, \eqref{final_cond2} and the condition given by Lemma \ref{harman_s2}, we get the following linear optimization system: setting $x=\frac1k$ and let $w$ be the exponent of $\eta$ we would like to optimize,
\begin{align*}
\begin{cases}
&x\le1;\ w\ge0;\ u\le\frac{1}{14}\\
&-w\ge\frac12-\frac{x}{2}-u\\
&-w\ge\frac14-\frac{x}{2}+2u.\\
\end{cases}
\end{align*}

Solving the system, it turns out that $u=\frac{1}{14}$ (and consequently $X=q^{7/3}$) are the optimal values; unfortunately, condition \eqref{final_cond2} does not affect linear optimization for values of $u\le 1/14$. Then, the maximum $k$-range is $\left(1,\frac{7}{6}\right)$ and
\begin{align*}
\eta=\left(\max_j p_j\right)^{-\frac{7-6k}{14k}+\varepsilon}.
\end{align*}

We thank the anonymous referees for an extremely careful reading of a previous version of this paper and the fruitful suggestions.


\begin{thebibliography}{99}
%
%

\bibitem{brudern-cook-perelli}
Br\"udern, J., Cook, R., Perelli, A., The values of binary linear forms
  at prime arguments. Proc of Sieve Methods, Exponential Sums and Their
  Application in Number Theory, Cambridge: Cambridge Univ Press, 87--100 (1997).

\bibitem{baker-harman}
Baker, R., Harman, G., Diophantine approximation by prime numbers.
  Journal of the London Mathematical Society 2~(2), 201--215 (1982).

\bibitem{cook2001general}
Cook, R., The value of additive forms at prime arguments. Journal de
  Th{\'e}orie des Nombres de Bordeaux 13~(1), 77--91 (2001).

\bibitem{cook-fox}
Cook, R., Fox, A., The values of ternary quadratic forms at prime
  arguments. Mathematika 48~(1-2), 137--149 (2001).

\bibitem{cook-harman}
Cook, R., Harman, G., The values of additive forms at prime arguments.
  The Rocky Mountain Journal of Mathematics 36~(4), 1153--1164 (2006).

\bibitem{davenport1946indefinite}
Davenport, H., Heilbronn, H., On indefinite quadratic forms in five
  variables. Journal of the London Mathematical Society 1~(3), 185--193 (1946).
  
\bibitem{Gambini}
Gambini, A., Diophantine approximation with one prime, two squares of primes and one $k$-th power of a prime. Open Mathematics, 19(1) 373--387 (2021).
  
\bibitem{GLZ}
Gambini, A., Languasco, A., Zaccagnini, A., A diophantine approximation problem with two primes and one $ k $-power of a prime. Journal of Number Theory 188, 210--228 (2018).

\bibitem{ge-li}
Ge, W., Li, W., One {D}iophantine inequality with unlike powers of prime
  variables. Journal of Inequalities and Applications ~(1), 1--8 (2016).  
  
\bibitem{gao-liu}
  Gao, Gaiyun, and Zhixin Liu. "Results of Diophantine approximation by unlike powers of primes." Frontiers of Mathematics in China 13 (2018): 797-808.
  

\bibitem{harman1991general}
Harman, G., Diophantine approximation by prime numbers. Journal of the
  London Mathematical Society 2~(2), 218--226 (1991).

\bibitem{harman2004ternary}
Harman, G., The values of ternary quadratic forms at prime arguments.
  Mathematika 51~(1-2), 83--96 (2004).
  
\bibitem{harman-kumchev}
Harman, G., Kumchev, A.V., On sums of squares of primes. Mathematical Proceedings of the Cambridge Philosophical Society, Cambridge Univ Press 140, 1--13 (2006).


\bibitem{languasco-zaccagnini-quaternary}
Languasco, A., Zaccagnini, A., A {D}iophantine problem with a prime and
  three squares of primes. Journal of Number Theory 132~(12), 3016--3028 (2012).

\bibitem{languasco-zaccagnini-ternary}
Languasco, A., Zaccagnini, A., On a ternary {D}iophantine problem with
  mixed powers of primes. Acta Arithmetica 159~(4), 345--362 (2013).

\bibitem{Languasco-Zaccagnini2016}
Languasco, A., Zaccagnini, A., A {D}iophantine problem with prime
  variables. V. Kumar Murty, D. S. Ramana, and R. Thangadurai, editors, Highly Composite: Papers in Number
Theory, Proceedings of the International Meeting on Number Theory, celebrating the 60th Birthday of Professor R. Balasubramanian (Allahabad, 2011), volume 23, pages 157--168. RMS-Lecture Notes Series (2016).

\bibitem{li-wang}
Li, W., Wang, T., Diophantine approximation with one prime and three squares of primes. The Ramanujan Journal 25~(3), 343--357 (2011).
  
\bibitem{liu2013diophantine}
Liu, Z., Sun, H., Diophantine approximation with one prime and three
  squares of primes. The Ramanujan Journal 30~(3), 327--340 (2013).

\bibitem{matomaki}
Matom{\"a}ki, K., Diophantine approximation by primes. Glasgow
  Mathematical Journal 52~(01), 87--106 (2010).

\bibitem{mu2016diophantine}
Mu, Q., Diophantine approximation with four squares and one $k$-th power
  of primes. The Ramanujan Journal 39~(3), 481--496 (2016).
  
\bibitem{muzhuli}
  Mu, Quanwu, Minhui Zhu, and Ping Li. A Diophantine inequality with four squares and one k th power of primes. Czechoslovak Mathematical Journal 69 (2019): 353-363.
  
\bibitem{rieger}
Rieger, G., \"Uber die Summe aus einem Quadrat und einem Primzahlquadrat. J. reine angew. Math. 231, 89--100 (1968).

\bibitem{vaughan1974diophantineI}
Vaughan, R., Diophantine approximation by prime numbers,
  {I}. Proceedings of the London Mathematical Society 3~(2), 373--384 (1974).

\bibitem{vaughan1974diophantineII}
Vaughan, R., Diophantine approximation by prime numbers,
  {II}. Proceedings of the London Mathematical Society 3~(3), 385--401 (1974).


\bibitem{wang-yao}
Wang Y., Yao W., Diophantine approximation with one prime and three squares of primes. Journal of Number Theory 180, 234--250 (2017).

\bibitem{watson}
Watson, G., On indefinite quadratic forms in five variables. Proceedings
  of the London Mathematical Society 3~(1), 170--181 (1953).

\end{thebibliography}
\end{document}